\documentclass[12pt]{amsart}
\usepackage{amsmath,amssymb, euscript, graphicx, hyperref}

\theoremstyle{definition}

\newcommand{\blackboard}[1]{\ensuremath{\mathbb{#1}}}

\newcommand{\N}{\blackboard{N}}

\newcommand{\Z}{\blackboard{Z}}

\newcommand{\R}{\blackboard{R}}

\begin{document}

\author{Azer Akhmedov}
\address{Azer Akhmedov,
Department of Mathematics,
North Dakota State University,
PO Box 6050,
Fargo, ND, 58108-6050}
\email{azer.akhmedov@ndsu.edu}

\title{Big tiles in hyperbolic groups}

\begin{abstract}  We prove that if $\Gamma $ is a word hyperbolic group and $K$ is a finite subset of $\Gamma $, then $\Gamma $ admits a tile containing $K$.
\end{abstract}

\maketitle


\
  \medskip
 
 A {\em tile} in a countable group $\Gamma $  is any finite subset $F$ of $\Gamma $ of cardinality at least two such that with non-overlapping left shifts of $F$ one can cover the whole group $\Gamma $. If $F$ is a tile of $\Gamma $,  then $\Gamma  $ admits a partitioning  $\Gamma = \displaystyle \mathop{\sqcup }_{
g\in C }gF$ for some subset $C\subset \Gamma $ called the {\em center set} of the {\em  tiling} $\tau(\Gamma , F) = \displaystyle \mathop{\sqcup }_{
g\in C  }gF$. A tile may have more than one (in fact, uncountably many) tiling. 

 \medskip 

 The following property of groups is studied  in \cite{Ch} and \cite{O-W}:

\medskip 

$(P)$ : Suppose that $\Gamma $ is a countable group. We say that $\Gamma $ has property $(P)$ if and only if  for all finite
subsets $K$ of $\Gamma $, there is a tile of $\Gamma  $ containing $K$.

\medskip 

{\bf Example 1.} A 2-element subset $\{x, y\}$ is a tile of a countable group $\Gamma $ if and only if the element $x^{-1}y$ is either non-torsion or has even order.

\medskip 

{\bf Example 2.} The set $\{0, 1, 3\}$ is not a tile of $\Z $. More generally, the set $\{0, 1, x\}$
is a tile if and only if $x \equiv 2 (mod \ 3)$. On the other hand, any finite arithmetic progression in $\Z $ is a tile of $\Z $.

\medskip 

{\bf Example 3.} Let $p$ be a prime number, $A = \{a_0, a_1, a_2, \dots , a_{p-1}\} \subset \Z$ be a
finite subset where $a_0 < a_1 < \dots  < a_{p-1}$ and $A(z) = \displaystyle \mathop{\sum }_{i=0}^{p-1} z^{b_i}$
where $b_i = a_i - a_0$ for all $i \in  \{0, 1, \dots , p-1\}$. Then $A$ is a tile if and only if there
exists $k\in \N $ such that the cyclotomic polynomial $\Phi _{p^k} (z)$ divides $A(z)$, see \cite{N}. Newman \cite{N} has
determined all tiles $A$ of $\Z $ where $|A|$ is a prime power. In general, however, the
problem remains open. A positive solution to the following conjecture would be a nice characterization of tiles in $\Z $ ($\Z $-tiles).

 \medskip

 {\bf Conjecture of Coven-Meyerowitz, \cite{C-M}, 1999} : Let $A \subset \Z$ be a finite
subset, $R_A = \{d \in \N | \Phi _d  \ \mathrm{divides} \  A(x)\}$, and $S_A = \{p^{\alpha } \in R_A\}$ - the set of prime
powers of $R_A$. Then $A$ is a $\Z $-tile if and only if the following conditions are satisfied:

(T1) : $B(1) = \Pi _{p^{\alpha }\in S_A} p$,

(T2) : if $x_1, x_2, \dots , x_n \in S_A$ then $x_1x_2 \dots x_n \in R_A$.

\medskip 

 It is known that the conditions  $(T_1)$  and  $ (T_2)$ together imply that   $A$  is  a $\Z $-tile. Moreover, if $A$ is a $\Z $-tile, then condition $(T_1)$ holds.
It is not known whether or not if $A$  being  a $ \Z$-tile implies the condition $(T_2)$.

\medskip 

{\bf Example 4.}  Balls in the Cayley graphs of groups may easily fail to be tiles. For the group $\Z^d =\langle a_1, \dots , a_d  \ | \ a_ia_j=a_ja_i, 1\leq i < j \leq d \rangle $, for any $n\geq 1$, the ball $B_n$ of radius $n$ with respect to the generating set $\{a_1, \dots , a_d\}$ as well as the set $[n]^d : = \{a_1^{i_1}a_2^{i_2}\dots a_d^{i_d} \ | \  1 \leq i_1, \dots , i_d \leq n\}$ are natural tiles. In the recent work of Benjamini-Kozma-Tzalik \cite{B-K-T}, the number $t_{n,d}$ of subsets $[n]^d$ which tile $\Z^d$ is estimated as $$(\sqrt[3]{3})^{n^d-o(n^d)}\leq t_{n,d}\leq (\sqrt[3]{3})^{n^d+o(n^d)}.$$
\medskip 

For general $\Z^n$, we would like to mention the following

\medskip 

{\bf Fuglede Conjecture \cite{F}, 1974} : If $A \subseteq \R^n$
is a measurable subset then $A$ is a tile of $\R ^n$
if and only if $A$ is spectral, i.e. for some set (spectrum) $\Lambda  \subseteq \R ^n$, the space $L^2(A)$ has an orthogonal basis $\{e^{2\pi i\lambda x} \ | \  \lambda\in \Lambda \}$. 

\medskip 

 The conjecture has been disproved by T.Tao \cite{T} for $n\geq 5$; and \cite{M}] and \cite{F-M-M} extended the result to the cases of $n = 4$ and $n=3$ respectively. For $n = 1$ and for $n = 2$, it still remains open.

\medskip 

{\bf Example 5}. Let $r \in \N$ and $B_r$ be a ball of radius $r$ in the Cayley graph of the free group $\mathbb{F}_k = \langle a_1, \dots , a_k \rangle $ with respect to standard generating set. Then the
sets $B_1$ and $B_1\backslash \{1\}$ are tiles of $\mathbb{F}_k$. On the other hand, it is not difficult to see that $B_r$ is a tile for every $r\geq 2$ as well, whereas  $B_r\backslash \{1\}$ is not a tile for $r\geq 2$.

\medskip 

{\bf Example 6.} If $\mathbb{F}_k = \langle  a_1, \dots , a_k \rangle $ and $k\geq 2$, then a sphere $S_r$  in the Cayley graph with respect.to the standard generating set is a tile if and only if $r = 1$. 

\medskip 

{\bf Example 7.} Any connected set in the Cayley graph of $\mathbb{F}_k = \langle  a_1, \dots , a_k \rangle $ with respect to the standard generating set is a tile \cite{A-F}. 

\bigskip 

In \cite{Ch} the following facts are observed and proved:

\medskip 

(i) If $G_1$ is a normal subgroup of $G$ and if both $G_1$ and $G/G_1$ have property $(P)$, then $G$ has property $(P)$;

\medskip 

(ii) If $G_n, n\geq 1$ have property $(P)$ where $G_n\leq G_{n+1}$ for all $n\geq 1$, then $\displaystyle \mathop{\cup }_{n\geq 1}G_n$ satisfies property $(P)$;

\medskip 

(iii) Abelian groups have property $(P)$.

\medskip 

 Since a set of representatives of all cosets with respect to a given finite index subgroup form a tile, we can add the following statement:    

\medskip 

(iv) Residually finite groups have property $(P)$.

\medskip 

As a corollary of (i) and (iii), one can state that (as it is done in \cite{Ch}) solvable groups have property $(P)$. Let us recall that not all (finitely generated) solvable groups are residually finite. As an immediate  corollary of (i), (ii) and (iii),  it is also observed in \cite{Ch}] that elementary amenable groups have property $(P)$. Properties (i)-(iii) also allow us to observe that the property $(P)$ is preserved under taking direct sums, direct products, wreath products, etc. Let us also observe that it is preserved under taking the free product (this is shown also in \cite{Ch} by a different argument):

\medskip 

 (v) If $G, H$ have property $(P)$, then so does their free product $G\ast H$.  

\medskip 

 Indeed, we have a short exact sequence $$1\rightarrow [G,H]\rightarrow G\ast H \rightarrow G\times H\rightarrow 1$$
where $[G,H]$ is the normal subgroup generated by the set $\{[x,y] : x\in G, y\in H\}$. But the subgroup $[G,H]$ is free (cf. \cite{S}, Prop.4, page 6) and free groups have property $(P)$. (It suffices to verify property $(P)$ for finitely generated subgroups and finitely generated free groups satisfy $(P)$ since, with respect to the standard generating set, any connected set is a tile \cite{A-F} or alternatively, any ball is a tile; see Example 5.) 
 \medskip 

 The following question, asked in \cite{Ch} and \cite{O-W}, is still open: {\em Does every countable group have property $(P)$?}

 \medskip

 The purpose of this paper is to answer this question positively for word hyperbolic groups. Let us remind a reader that a finitely generated group is word hyperbolic if its Cayley graph with respect to some (consequently, with respect to any) finite generating set is a hyperbolic space in the sense of Gromov (see \cite{D}, \cite{K-B} or \cite{G-H} for this and other basic notions of the theory of hyperbolic groups and hyperbolic spaces in general). Let us recall that for a geodesic metric space $X$, we say $X$ is $\delta $-{\em hyperbolic} for some non-negative $\delta \geq 0$ if all triangles are $\delta $-thin. This means that for any three points $p_i\in X, i\in \Z/3\Z$ and geodesics $\alpha _i$ (viewed as subsets of $X$) connecting $p_{i+1}$ and $p_{i-1}$, for all $i\in \Z/3\Z$, any point of $\alpha _i$ has distance  at most $\delta $ from the union $\alpha _{i-1}\cup \alpha _{i+1}$. The metric space $X$ is called {\em hyperbolic} if it is $\delta $-hyperbolic for some $\delta \geq 0$.

\medskip 

 We will also use another well-known description of thinness of triangles. For this purpose, for positive real numbers $p,q,r$, let $T_{p,q,r}$ be a tree with one vertex of degree 3 (called the {\em central vertex}) and three vertices of degree 1 (called {\em leaves}) endowed with the tree metric where the edges have length $p, q$ and $r$. If the triangle $\Delta = \alpha _1\cup \alpha _2\cup \alpha _3$ (as in the previous paragraph), consisting of geodesics $\alpha _1, \alpha _2, \alpha _3$ is $\delta $-thin, then there exist points $m_i\in \alpha _i, i\in   \Z/3\Z$ such that  $$\max \{d(m_i, m_{i+1}) :  i\in \Z/3\Z\}\leq 4\delta , d(p_i, m_{i+1}) = d(p_i, m_{i-1}) : = l_i,$$ and there exists a continuous map $\phi :\Delta \to T_{l_1, l_2, l_3}$ where  $T_{l_1, l_2, l_3}$ is consisting of geodesics $\eta _0, \eta _1, \eta _2$ of length $l_0, l_1, l_2$ respectively connecting the central point $a$ with the leaves $v_0, v_1, v_2$ respectively such that  the restriction of $\phi $ to each side of $\Delta $ is an isometry; in addition,  $$\phi ^{-1}(a) =  \{m_1, m_2, m_3\} \ \mathrm{and} \ \phi ^{-1}(v_i) = \{p_i\}, i\in \Z/3\Z $$ and for all $i\in \Z /3\Z, x\in \eta _i\backslash \{a, v_i\}$, we have $\phi ^{-1}(x) = \{y, z\}$ where the point $y$  belongs to the part of geodesic $\alpha _{i+1}$ connecting $p_i$ and  $m_{i+1}$, the point $z$ belongs to the part of geodesic $\alpha _{i-1}$ connecting $p_i$ and  $m_{i-1}$, and $d(y, z) \leq 4\delta .$   A geodesic triangle with vertices at $m_0, m_1, m_2$ will be called a {\em core} of the triangle $\Delta $. Such a core will  lie in a ball of radius $4\delta $.

\bigskip 

We now state our key result.

\medskip 

{\bf Proposition 1.} Let $\Gamma $ be a non-elementary word hyperbolic group without a non-trivial finite normal subgroup. Then for any finite subset $K \subset \Gamma $ there exists $\xi \in \Gamma $ such that $K \cup  \{\xi \}$ is a tile of $\Gamma $.

\medskip 

 For the proof of this proposition, we need to recall some facts about word hyperbolic groups. Let $\Gamma $ be a word hyperbolic group with a fixed finite generating set $S$ where $S = S^{-1}$ and $1\notin S$. Then the Cayley graph with respect to this generating set is $\delta $-hyperbolic for some $\delta > 0$. For any $g\in \Gamma $, we will write $|g|$ for the shortest length of a path representing $g$, and $d(g_1, g_2) = |g_1^{-1}g_2|$ for the corresponding left-invariant Cayley metric. 

\medskip 

  $\Gamma $ acts on $\Gamma $ by left translations which are isometries. In addition, $\Gamma $ admits a boundary $\partial \Gamma $ which is a compact metric space upon which $\Gamma $ acts by homeomorphisms. $\Gamma \cup \partial \Gamma $ is also compact and can be viewed as a compactification of $\Gamma $.  

\medskip

  Non-identity isometries of a hyperbolic space $X$ can be classified as {\em hyperbolic} (elements with two fixed points on the boundary $\partial X$ one of which is attractive and another one repelling), {\em parabolic} (elements with one fixed point on the boundary $\partial X$ where attractive and repelling points coincide) and {\em elliptic} (all other non-identity isometries; equivalently, an elliptic isometry is a non-identity isometry such that the orbit of any point of $X$ under the iterations of it is bounded). In the restricted context when $X$ is the Cayley graph of $\Gamma $, a non-torsion element acts as a hyperbolic element and torsion elements are elliptic (so we have no parabolic elements). A non-torsion element $\gamma \in \Gamma $ has two fixed points $P_a(\gamma )$ and  $P_r(\gamma )$ on $\partial \Gamma $; $P_a(\gamma )$ ($P_r(\gamma )$) is an attractive (repelling) point, i.e. for all open neighborhoods  $U, V$ of $P_a(\gamma )$ and $P_r(\gamma )$ respectively, there exists $N\geq 1$ such that for all $n > N$, $\gamma ^n(\partial \Gamma \backslash V)\subseteq U$ ($\gamma ^{-n}(\partial \Gamma \backslash U)\subseteq V$). A torsion element may have no fixed point or many (finitely many or infinitely many) fixed points; in particular, an elliptic element may fix all points of $\partial \Gamma $. We will collect several well known facts in the following

   \medskip 

   {\bf Proposition 2.} Let $\Gamma $ be a non-elementary word hyperbolic group. Then 

   a) $\partial \Gamma $ is infinite;

   b) the sets $\{P_a(\gamma ) : \gamma \ \mathrm{is \ a \ non}$-$\mathrm{torsion \ element}\}$ and $\{P_r(\gamma ) : \gamma \ \mathrm{is \ a \ non}$-$\mathrm{torsion \ element}\}$ are dense in $\partial \Gamma $;

   c) the set $\{(P_a(\gamma ), P_r(\gamma )) : \gamma \ \mathrm{is \ a \ non}$-$\mathrm{torsion \ element}\}$  is dense in $\partial \Gamma \times \partial \Gamma$;

   d) The set $\{\gamma \in \Gamma : \gamma \ \mathrm{fixes} \ \partial \Gamma \ \mathrm{pointwise}\}$ forms a finite normal subgroup. \ $\square $

   \medskip 

   Using the above proposition, we will prove the following 

   \medskip 

   {\bf Lemma 1.}  Let $\Gamma $ be a non-elementary word hyperbolic group with a fixed generating set and without a non-trivial finite normal subgroup, and $B$ be a finite subset of $\Gamma \backslash \{1\}$. Then for all $N\geq 1$ there exists a non-torsion element $\xi \in \Gamma $ such that   $|\xi ^{mi}b\xi ^{mj}| > N$ for all $m\geq 1, i,j\in \{-1,1\}$ and $b\in B$.

   \medskip 

   {\bf Proof.} Let $A = B\cup B^{-1}\cup \{x\in \Gamma : 0 < |x| \leq N\}$ and $A = \{a_1, \dots , a_n\}$. Since $\Gamma $ has no finite normal subgroup, for all $1\leq s\leq n$, there exists $p_s\in \partial \Gamma $ such that $a_s(p_s)\neq p_s$. 

\medskip 

    Let $P = \{p_1, \dots , p_n\}$. By Proposition 2, there exists a hyperbolic $\eta \in \Gamma \backslash \{1\}$ such that $\{a(P) : a\in A\cup \{1\}\}\cap Fix(\eta ) = \emptyset $.  Since $\eta $ is non-torsion, we have $|\eta ^n|\to \infty $ as $n\to \infty $. If the claim of the lemma does not hold for any power $\xi = \eta ^l$ of $\eta $, then we have $|\eta ^{mi}b\eta ^{mj}| \leq N$ for some $b\in B$ and for infinitely many $m\geq 1$. Then, necessarily, for all $M\geq 1$, there exist $b\in A, i,j,k \in \{-1,1\}$ and $p, q\geq M$ such that $\eta ^{ip}b^k\eta ^{jq} = b $. Then, since $A$ is finite, there exist $b\in A, i,j,k \in \{-1,1\}$ such that for all $M\geq 1$ there exist $p, q\geq M$ such that  $\eta ^{ip}b^k\eta ^{jq} = b \ (\dagger )$.

     \medskip 

     Let $b = a_r$  for some  $1\leq r\leq n$. By our choice,  $\{p_r, b(p_r), b^{-1}(p_r)\}\cap Fix(\eta ) = \emptyset $.  But for any open set $U\supseteq Fix(\eta )$, if $p$ and $q$ are sufficiently big, then $\eta ^{ip}b^k\eta ^{jq}\in U$ thus using $(\dagger )$ we obtain a contradiction by taking $U$ containing  $Fix(\eta )$ and disjoint from $ \{p_r, b(p_r), b^{-1}(p_r)\}$. \ $\square $
     
    \bigskip

{\bf Proof of Proposition 1.} Let us fix a finite generating set of $\Gamma $. Then the Cayley graph with respect to this generating set is $\delta $-hyperbolic for some $\delta > 0$. For any $g\in \Gamma $,  $|g|$ will denote the length of a shortest path representing $g$, and $d(g_1, g_2) = |g_1^{-1}g_2|$ for the corresponding left-invariant Cayley metric. For any $g_1, g_2\in \Gamma $, $[g_1, g_2]$ will denote the set of all geodesics connecting $g_1$ and $g_2$.  

\medskip 

 If $|K| = 1$, then assuming $K = \{u\}$ we can choose $\xi = ux$ where $x$ is any non-torsion element. (Let us recall a well known fact that a word hyperbolic group has finitely many conjugacy classes of torsions and any non-finite word hyperbolic group contains a non-torsion element.) So, let us assume that $|K|\geq 2$ and $\{1, v\} \subseteq K$. Let $$r = 4\max\{\max \{|g| : g\in K\cup K^{-1}\}, [\delta +1]\}+1.$$

Because of $\delta $-thinness of triangles, there exists $\xi  \in \Gamma $ such that the set $$M_{\xi ,r} = \{x \in  \Gamma  : |x\xi | < |x| + r \ \mathrm{and} \  |x\xi ^{-1}| < |x| + r\}$$ is empty and the set $$S_{\xi ,r} = \{x \in \Gamma  : |x\xi | < |x| + r \ \mathrm{or} \ 
|x\xi ^{-1}| < |x| + r\}$$ is $4r$-separated, i.e. for any $y, z \in  S_{\xi ,r}, d(y,z) > 4r$.

 \medskip 

  Indeed, by Lemma 1, we can choose a non-torsion $\xi $ with $|\xi | > 100r$ such that for all $m\geq 1, i,j\in \{-1,1\}$, we have $|\xi ^{im}b\xi ^{jm}| > 100r$ for all $b$ with $1\leq |b| \leq 4r$. By replacing $\xi $ with $\xi ^m$ for sufficiently big $m$ if necessary, we may assume that for all $b$ with $1\leq |b| \leq 4r$ and for all $n\geq 1, i, j\in \{-1,1\}$, $|\xi ^ib\xi ^j| > (1.9) |\xi |$. In addition, again by replacing $\xi $ with $\xi ^m$ for sufficiently big $m$ if necessary, we can also claim the inequality $|\xi ^2| > 1.9|\xi |$. Indeed, since $\xi $ is not a torsion, there exist  natural numbers $l_0, L$ such that $|\xi ^{l_0}|\geq 10(\delta +1) \ (\dagger _1)$ and for all natural $m\geq L$,  $|\xi ^{m}|\geq 200|\xi ^{l_0}| \ (\dagger _2)$. For all $m\geq L$ consider the triangle $\Delta _m$ with vertices at $1, \xi ^m, \xi ^{2m}$ and sides $\alpha \in [1, \xi ^m], \beta \in [\xi ^m, \xi^{2m}]$ and $\gamma \in [1, \xi^{2m}]$, and let $C(\Delta _m)$ be a core of $\Delta _m$. If $a, b$ are the vertices of $C(\Delta _m)$ on the sides $\alpha , \beta $ respectively, then because of $(\dagger _1)$ and $\delta $-thinness of triangles, we have $\max \{d(a, \xi ^m), d(b, \xi ^m)\} < |\xi ^{l_0}|$. Then, because of $(\dagger _2)$, we have $|\xi ^{2m}| > 1.9|\xi ^{m}|$. Thus we established that for all $m\geq L$, $|\xi ^{2m}| > 1.9|\xi ^{m}|$. 

  \medskip 

  Now, to see that $M_{\xi ,r} = \emptyset $, let us consider a geodesic triangle $\Delta $ with vertices at $1, x\xi , x\xi ^{-1}$ with sides $\alpha \in [1, x\xi] , \beta \in [1,x\xi ^{-1}], , \gamma \in [x\xi , x\xi ^{-1}]$ where $x$ is any element of $M_{\xi ,r}$. Then, by inequality $|\xi ^2| > (1.9)|\xi |$ and by $\delta $-thinness of a triangle with vertices at $x, x\xi , x\xi ^{-1}$, there exists a point $p$ on $\gamma $ such that $$d(p,x) < 0.2|\xi | \ \mathrm{and} \ \min\{d(p, x\xi), d(p, x\xi ^{-1})\} > 0.8|\xi | \ \ \ (1)$$ 

   Indeed, let $\gamma _1\in [x\xi , x], \gamma _2\in [x\xi ^{-1},x]$ and $\Delta $ be a geodesic triangle with vertices $x , x\xi ^{-1}, x\xi $ and sides $\gamma , \gamma _1, \gamma _2$. Let also $\Delta _0$ be a core of $\Delta $ with vertices $g\in \gamma ,g_1\in \gamma _1, g_2\in \gamma _2 $. We will show that the point $p$ can be taken as $g$ to satisfy inequalities (1). 

   \medskip 

  \begin{center}
      \includegraphics[width=12cm, height=8cm]{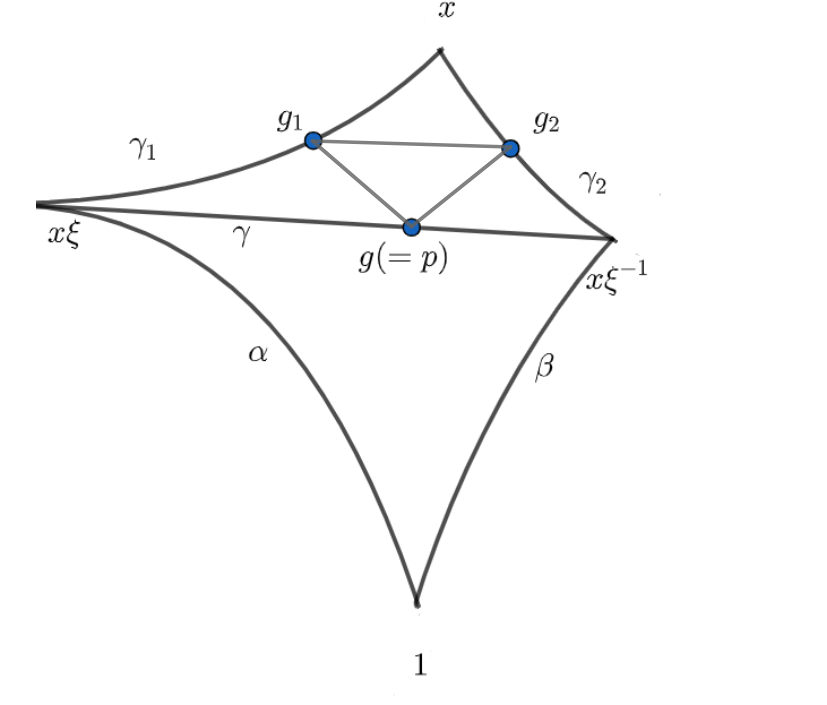}
  \end{center} 

   \medskip

    We have $$1.9|\xi |\leq |\xi ^2| = d(\xi, \xi ^{-1}) =  d(x\xi, x\xi ^{-1})\leq d(x\xi , g_1)+d(g_1,g_2)+d(g_2,x\xi ^{-1})\leq $$ \ $$\leq d(x\xi , g_1)+4\delta + d(g_2,x\xi ^{-1}) = (|\xi |-d(g_1,x))+4\delta + (|\xi |-d(g_2,x)) = $$ \ $$ 2|\xi |-2d(g_1,x)+4\delta .$$

\medskip 

     So we obtain that $d(g_1,x)\leq 0.1|\xi | + 2\delta  \  (2).$

      By (2), we also have $$d(g,x\xi ) = d(g_1,x\xi ) = d(x,x\xi )-d(g_1, x) = |\xi | - d(g_1, x) \geq |\xi | - 0.1|\xi | - 2\delta > 0.8|\xi |$$ thus $d(g,x\xi ) > 0.8|\xi | \ (3)$. Similarly, we obtain that $d(g,x\xi ^{-1}) > 0.8|\xi | \ (4)$.  

      \medskip 

       On the other hand, using inequality (2) again, we have $$d(g,x) \leq d(g,g_1) + d(g_1,x) \leq 4\delta + d(g_1, x) \leq 0.1|\xi | + 6\delta < 0.2|\xi | \ \ \  (5).$$

       \medskip 

       Now, we let $p = g$; then the inequalities (3), (4) and (5) yield the inequalities (1).

   \medskip

  But by $\delta $-thinness of triangles, $p$ must be at most $\delta $-apart either from $\alpha $ or from $\beta $; without loss of generality, we may assume that $d(p, \alpha ) \leq \delta $ and $q$ is a point on $\alpha $ with $d(p,q)\leq \delta $. Then applying the triangle inequality, we have $$d(1,x\xi ) = d(1,q)+d(q,x\xi ) \geq (d(1,p)-d(p,q))+(d(p,x\xi )-d(p,q)) \geq $$ \ $$d(1,p)+d(p,x\xi )-2\delta \ (6).$$

       \medskip 

       Now, using (6), we obtain that $$|x|+r \geq |x\xi | = d(1,x\xi ) \geq d(1,p)+d(p, x\xi ) - 2\delta \geq (|x| - 0.2|\xi |) + 0.8|\xi | - 2\delta $$ \ $$ = |x| + 0.6|\xi | - 2\delta $$  
which implies that $r+2\delta \geq 0.6|\xi |$. Contradiction; thus $M_{\xi ,r} = \emptyset $.  

\medskip 

 Establishing the $4r$-separatedness of $S_{\xi ,r}$ is also similar. Indeed, let $y, z\in S_{\xi ,r}$ with $d(y,z) \leq 4r$ such that $|y\xi ^i| < |y| + r$ and $|z\xi ^j| < |z| + r$ for some $i, j\in \{-1,1\}$. We let $\alpha \in [1, y\xi ^i], \beta \in [1, z\xi ^j], \gamma \in [y\xi ^i, z\xi ^j]$. Then, again, by the inequality $|\xi ^ib\xi ^j| > (1.9) |\xi |$ and by $\delta $-thinness of triangles, there exists $p$ on $\gamma $ such that $$\max \{d(p,y), d(p,z)\} < 0.2|\xi | \ \mathrm{and} \ \min\{d(p, y\xi ^i), d(p, z\xi ^j)\} > 0.8|\xi | \ (7)$$ 

   \medskip

   Indeed, let $y_1 = y\xi ^i, z_1 = z\xi ^j$. Then $y_1^{-1}z_1 = (y_1^{-1}y)(y^{-1}z)(z^{-1}z_1) = \xi ^{-i}(y^{-1}z)\xi ^j$ hence by the inequality $|\xi ^{-i}b\xi ^j| > 1.9|\xi |$, we have $d(y_1,z_1) > 1.9|\xi |$. Then, defining the triangle $\Delta $ and the points $g, g_1, g_2$ as before, we can write $$1.9|\xi | < d(y_1, z_1)\leq d(y_1, g_1)+d(g_1,g_2)+d(g_2,z_1) \leq d(y_1, g_1)+4\delta +d(g_2,z_1) = $$ \ $$(d(y_1,x)-d(x,g_1))+4\delta +(d(z_1,x)-d(x,g_2) = d(y_1,x)+d(z_1,x) +4\delta -2d(x,g_1) \leq $$ \ $$(|\xi |+4r)+(|\xi |+4r)+4\delta -2d(x,g_1).$$

\medskip 

         So $d(x,g_1) < 0.1|\xi |+4r+2\delta \ (8)$. Then $$d(y,g) \leq d(y,x)+d(x,g_1)+d(g_1,g)\leq 4r + (0.1|\xi |+4r+2\delta ) + 4\delta  = $$ \ $$8r+0.1|\xi |+ 6\delta < 0.2|\xi |$$ thus $d(y,g) < 0.2|\xi | \ (9)$. Similarly, we obtain that $d(z,g) < 0.2|\xi | \ (10)$. On the other hand, $$d(g,y_1) \geq d(g_1,y_1) - d(g,g_1) \geq d(g_1,y_1)-4\delta = d(y_1,x)- d(g_1,x)-4\delta \geq $$ \ $$d(y_1,y) - d(y,x) - d(g_1,x)-4\delta \geq |\xi | - 4r  - d(g_1,x) - 4\delta $$ and from the inequality (8) we now obtain $d(g,y_1) > 0.8|\xi | \ (11)$. Similarly, we establish that $d(g,z_1) > 0.8|\xi | \ (12)$  

         \medskip 
         
          Now, by letting $p=g$, from the inequalities (9), (10), (11) and (12), we obtain the inequalities (7). 

   \begin{center}
      \includegraphics[width=12cm, height=8cm]{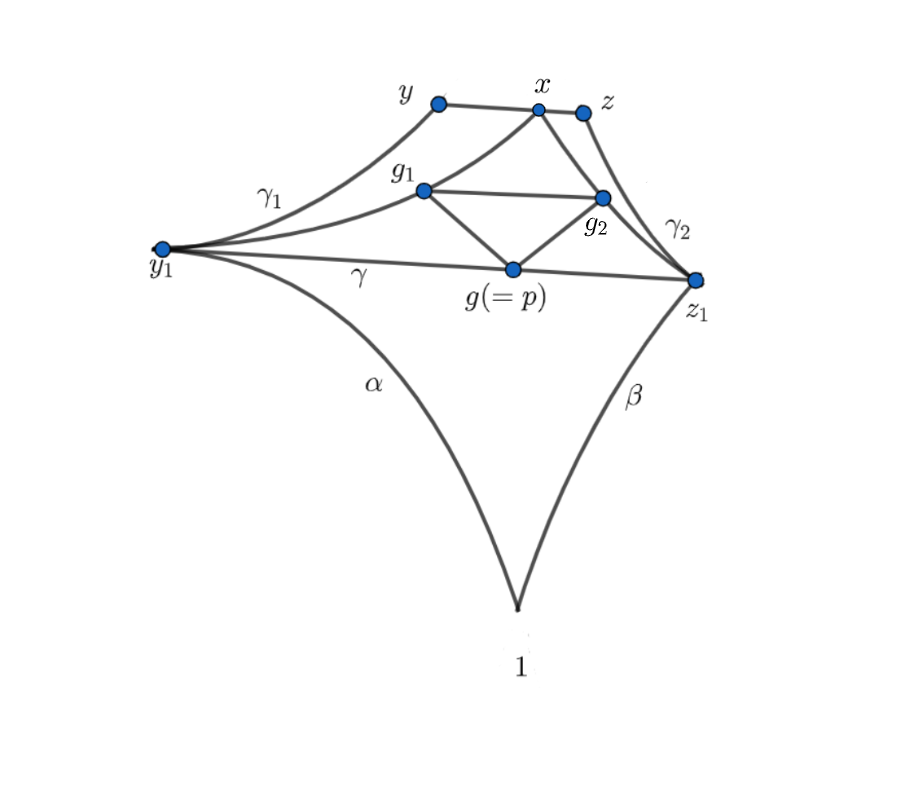}
  \end{center} 
 
 \medskip 
 
 Again, without loss of generality we may assume that $p$ is at most $\delta $ apart from $\alpha $ and using triangle inequalities as in (6) again, we have $d(1,y\xi ^i) \geq d(1,p)+d(p,y\xi ^i)-2\delta $ which yields $$|y|+r \geq |y\xi ^i| = d(1,y\xi ^i) \geq d(1,p)+d(p, y\xi ^i) - 2\delta \geq (|y| - 0.2|\xi |) + d(p, y\xi ^i) - 2\delta $$ \ $$\geq (|y| - 0.2|\xi |) + 0.8|\xi | - 2\delta  = |y| + 0.6|\xi | - 2\delta $$  
which implies that $r+2\delta \geq 0.6|\xi |$. Contradiction; thus we established the $4r$-separatedness of $S_{\xi , r}$. 

 \medskip 

 Now, let $C_{\xi ,r}= \{x\in S_{\xi ,r} : |x\xi ^{-1}| < |x| + r\}$. We first observe that the element $v^{-1}\xi $ is not a torsion. Indeed, assuming the opposite, let $(\xi^{-1}v)^q = 1$ for some $q\geq 1$ and $S$ be the set consisting of elements $1, \xi^{-1}, \xi^{-1}v, \xi^{-1}v\xi^{-1}, \dots , ( \xi^{-1}v)^{q-1}, ( \xi^{-1}v)^{q-1}\xi^{-1}$.

       \medskip 
       
       Let $\theta \in S$ such that $|\theta | = \displaystyle \mathop{\max}_{x\in S}|x|$. Then either $\theta = ( \xi^{-1}v)^i$ for some $1\leq i\leq q-1$ or 
       $\theta = ( \xi^{-1}v)^i\xi^{-1}$ for some $0\leq i\leq q-2$. In the former case, we have $|\theta \xi ^{-1}|\leq |\theta | < |\theta |+r$ and $|\theta v^{-1}\xi | \leq |\theta | < |\theta | + r$ thus $\theta , \theta v^{-1}\in S_{\xi , r}$, but this contradicts the fact that $S_{\xi , r}$ is $r$-separated. In the latter case, we have  $|\theta \xi |\leq |\theta | < |\theta |+r$ and $|\theta v\xi ^{-1}| \leq |\theta | < |\theta | + r$ thus $\theta , \theta v\in S_{\xi , r}$ which again contradicts $r$-separatedness of $S_{\xi , r}$.

 \medskip

 Thus we established that $v^{-1}\xi $ is not a torsion. Then $C_{\xi ,r}\subseteq S_{\xi ,r} $ and for all $s\in C_{\xi ,r}$ we have $sv^{-1}\xi \in C_{\xi ,r}, \{s, sv^{-1}\xi\} \subset sv^{-1}(K\cup \{\xi \})$, moreover, $C_{\xi ,r}\cap sv^{-1}(K\cup \{\xi \}) = \{s, sv^{-1}\xi\}$. (To see why $sv^{-1}\xi\in C_{\xi , r}$, notice that we have $v\neq 1$ and $|v| < r$. On the other hand, $s\in C_{\xi ,r} \subset S_{\xi, r}$. But $S_{\xi ,r}$ is $4r$-separated. So $sv^{-1}\notin S_{\xi , r}$. Then $|sv^{-1}\xi | \geq |sv^{-1}| + r.$   
    
    Then $$|(sv^{-1}\xi)\xi ^{-1}| = |sv^{-1}| \leq  |sv^{-1}\xi | - r <  |sv^{-1}\xi | + r.$$ 
 
     Hence $sv^{-1}\xi\in C_{\xi , r}$.) 
 
 Then we have a subset (and enumeration) $\{s_1, s_2, \dots \}\subset C_{\xi ,r}$ and a partition $C_{\xi ,r} = \displaystyle \mathop{\sqcup }_{n\geq 1}\{s_n, s_nv^{-1}\xi \}$ such that $|s_i| \leq |s_j|$ for all $i<j$.

 \medskip 
 
Let $A = \displaystyle \mathop{\cup }_{s\in C_{\xi ,r}} sv^{-1}(K \cup \{\xi \})$. By the definition of  $C_{\xi ,r}$, for any two distinct $i, j\geq 1$, we have $s_iv^{-1}(K \cup \{\xi \})\cap s_jv^{-1} (K \cup  \{\xi \}) = \emptyset $. We also have $A\supset C_{\xi ,r}$.

\medskip 

Let $x_1^{(1)}, x_2^{(1)}, \dots $  be all elements of the set $B = \Gamma \backslash A$ such that $|x_i^{(1)}| \leq |x_j^{(1)}|$ if $i < j$. We start an infinite process of covering the set $\Gamma \backslash A$ by the non-overlapping shifts of $K\cup \{\xi \}$, namely: 

\medskip 

 Let $B_1 = B, B_2 = B\backslash x_1^{(1)}\xi ^{-1}
(K \cup  \{\xi \})$. We denote the elements of $B_2$ by
$x_1^{(2)}, x_2^{(2)}, \dots $ so that $|x_i^{(2)}| \leq |x_j^{(2)}|$ if $i < j$, and let $$B_3 = B_2\backslash (x_1^{(1)}\xi ^{-1}(K \cup  \{\xi \})\cup x_1^{(2)}\xi ^{-1}(K \cup  \{\xi \})),$$ and continue the process to cover all elements of the set $B$, i.e. if the set $B_n$ and the elements $x_1^{(n)}, x_2^{(n)}, \dots $ are chosen, we set $$B_{n+1} = B_n\backslash (\displaystyle \mathop{\cup }_{1\leq i\leq n}x_1^{(i)}\xi ^{-1}
(K \cup  \{\xi \})$$ and denote
the elements of $B_{n+1}$ by $x_1^{(n+1)}, x_2^{(n+1)}, \dots $ so that $|x_i^{(n+1)}| \leq |x_j^{(n+1)}|$ if $i< j$. 

 \medskip 

 Let $x_1^{(i)}\xi ^{-1}(K \cup  \{\xi \})\cap x_1^{(j)}\xi ^{-1}(K \cup  \{\xi \}) \neq \emptyset $ for some distinct positive integers $i$ and $j$, and let $i<j$. Since $x_1^{(i)}\neq x_1^{(j)}$, we obtain one of the following three possibilities:

\medskip 

  (a) $x_1^{(i)}\in x_1^{(j)}\xi ^{-1}K$

 (b) $x_1^{(j)}\in x_1^{(i)}\xi ^{-1}K$

  (c) $x_1^{(i)}\xi ^{-1}K\cap x_1^{(j)}\xi ^{-1}K \neq \emptyset $

\medskip 

 In case (a), since $x_1^{(j)}\in B$, we obtain that for any $w\in K$, $$|x_1^{(j)}\xi ^{-1}w| \geq |x_1^{(j)}\xi ^{-1}| - |w| > |x_1^{(j)}| + r - r = |x_1^{(j)}| \geq |x_1^{(i)}|;$$ contradiction. 

 \medskip

Case (b) is also impossible, since by construction $(x_1^{(i)}\xi ^{-1}K)\cap B_{i+1} = \emptyset $ whereas $x_1^{(j)}\in B_j\subseteq B_{i+1}$. 

\medskip 

 Case (c) would imply that $0 < d(x_1^{(i)}\xi ^{-1}, x_1^{(j)}\xi ^{-1}) \leq 2r \ (\ast )$. On the other hand, since $x_1^{(i)}, x_1^{(j)}\in B$ and $C_{\xi, r}\subseteq A = \Gamma \backslash B$, we have that  $x_1^{(i)}, x_1^{(j)}\notin C_{\xi , r}$. Hence $|x_1^{(i)}\xi ^{-1}| \geq |x_1^{(i)}|+r$ and $|x_1^{(j)}\xi ^{-1}| \geq |x_1^{(j)}|+r$. Then $$|(x_1^{(i)}\xi ^{-1})\xi | \leq  |x_1^{(i)}\xi ^{-1}| - r < |x_1^{(i)}\xi ^{-1}| + r$$ and $$|(x_1^{(j)}\xi ^{-1})\xi | \leq  |x_1^{(j)}\xi ^{-1}| - r < |x_1^{(j)}\xi ^{-1}| + r$$ which implies $x_1^{(i)}\xi ^{-1}, x_1^{(j)}\xi ^{-1}\in S_{\xi, r}$, but since $S_{\xi, r}$ is $4r$-separated, this contradicts $(\ast )$. We again obtained a contradiction.

 \medskip 

 Thus we showed that the sets $x_1^{(n)}\xi ^{-1}(K \cup  \{\xi \}), n\geq 1$ are mutually disjoint and together cover $B$.  

 \medskip 
 
 On the other hand, if $x_1^{(i)}\xi ^{-1}(K\cup \{\xi \})\cap sv^{-1}(K\cup \{\xi \}) \neq \emptyset $ (i.e. $x_1^{(i)}\xi ^{-1}(K\cup \{\xi \})\cap A \neq \emptyset $) for some $i\in \N $ and $s\in \{s_1, s_2, \dots \}$, then by construction $x_1^{(i)}\notin sv^{-1}(K\cup \{\xi \})$ so either  $x_1^{(i)}\xi ^{-1}K\cap sv^{-1}K \neq \emptyset $ or $x_1^{(i)}\xi ^{-1}K\cap sv^{-1}\{\xi \} \neq \emptyset $. Either way, since  $\{s, sv^{-1}\xi\} \subset sv^{-1}(K\cup \{\xi \})$,  $M_{\xi , r} = \emptyset $ and $S_{\xi , r}$ is $4r$-separated, we obtain a contradiction. Thus the covering of $\Gamma = A\cup B$ is non-overlapping, and by construction every element of $\Gamma $ is covered. Thus, we obtain a tiling of $\Gamma $ by $K\cup \{\xi \}$ with a center set $\{s_iv^{-1} : i\geq 1\}\sqcup \{x_1^{(i)}\xi ^{-1} : i\geq 1\}$. \ $\square $

 \bigskip 

 {\bf Proposition 3.} If $\Gamma $ is a non-elementary word hyperbolic group, then $\Gamma $ satisfies property $(P)$.

{\bf Proof.} Indeed, let $H$ be a finite normal subgroup of $\Gamma $ (if no such subgroup  exists, then we are done by Proposition 1). Then, since $\Gamma $ does not contain an infinite direct limit of finite groups (let us recall a well known fact that any amenable subgroup of a word hyperbolic group is virtually cyclic),  there exists a finite subgroup $H_1 \triangleleft \Gamma $  such that $\Gamma _1 = \Gamma /H_1$ is a non-elementary word hyperbolic group with no non-trivial finite normal subgroup. By Proposition 1, $\Gamma _1$ satisfies property $(P)$. Since $H_1$ is finite, it also satisfies $(P)$. Then, by the observation (i) in \cite{Ch} as mentioned above, the group $\Gamma $, as an extension of $\Gamma _1$ by $H_1$, satisfies $(P)$.\ $\square $

\bigskip 

 We would like to remark that in a recent paper \cite{MM}, the technique and ideas of our paper  is adapted to a setting of acylindrically hyperbolic groups. As a result of this generalization, combined with other known results in literature, property $(P)$ is established for more examples of groups beyond the class of word hyperbolic groups. Most notably, it is shown that one-relator groups and two-dimensional Artin groups satisfy property $( P)$.

 \end{document}